\documentclass{amsart}

\usepackage{amssymb}

\theoremstyle{plain}
\newtheorem{theorem}{Theorem}
\theoremstyle{definition}
\newtheorem{definition}{Definition}

\newcommand\ZZ{{\mathbb{Z}}}
\newcommand\GL{{\mathsf{GL}}}
\newcommand\Aut{{\mathsf{Aut}}}
\renewcommand\SS{{\Sigma}}
\newcommand\Sd{{\Sigma_{\diamond}}}
\newcommand\Sdn{{(\Sigma_{\diamond})^n}}
\newcommand\ov{\overline}

\begin{document}

\title{Cayley automatic groups \\ are not Cayley biautomatic}

\author{Alexei Miasnikov}

\address{Department of Mathematical Sciences, Stevens Institute, Castle Point on Hudson, Hoboken, NJ 07030-5991, USA}
\email{amiasnikov@gmail.com}

\author{Zoran \v Suni\'c}

\address{Department of Mathematics, Texas A{\&}M University, MS-3368, College Station, TX 77843-3368, USA}
\email{sunic@math.tamu.edu}

\begin{abstract}
We show that there are Cayley automatic groups that are not Cayley biautomatic. In addition, we show that there are Cayley automatic groups with undecidable Conjugacy Problem and that the Isomorphism Problem is undecidable in the clas of Cayley automatic groups. 
\end{abstract}

\maketitle

\section{Introduction}

The notion of automatic groups, based on ideas of Thurston, Cannon, Gilman, Epstein and  Holt,  was introduced in~\cite{epstein-al:wp}. The initial motivation was to understand the fundamental groups of compact 3-manifolds and  make them tractable for computing. It was quickly realized that automatic group come short in dealing with manifolds of Nil and Sol types. This immediately triggered a search for suitable generalizations. In~\cite{bridson-g:indexed} Bridson and Gilman came up with a sufficiently powerful notion of automaticity (asynchronously automatic groups where regular languages are replaced with indexed languages) that covers all fundamental groups of compact 3-manifolds, but at the cost of losing all the nice algorithmic properties. 

Since 1990's  a lot of groups were  proved to be automatic (see the survey in \cite{kharlampovich-k-m:cayley}), but one frustration still lingers there. It turns out that many basic questions on automatic groups remain wide open despite a considerable effort by the group theoretic community. Three such basic problems ask if automatic groups are biautomatic, if they have a decidable Conjugacy Problem, and if the Isomorphism Problem is decidable within the class. The Cayley automatic groups, introduced in~\cite{kharlampovich-k-m:cayley} retain the basic algorithmic properties of the standard automatic groups (decidability of the Word Problem in quadratic time and decidability of the Conjugacy Problem in the biautomatic case) but form a much wider class of groups, which, in particular, contains  many  nilpotent and solvable groups, which are not automatic under the standard definition. From the algorithmic view-point this indicates that the new class gives a legitimate notion of automaticity. Another confirmation that Cayley automatic groups provide a robust  generalization of the standard automatic groups is given by the fact that the basic problems mentioned above can be tamed in this case.  Namely, we show that all three have a negative solution in the class of Cayley automatic groups. 

\begin{theorem}\label{t:biautomatic}
There are Cayley automatic groups that are not Cayley biautomatic. 
\end{theorem}

This answers a question raised implicitly at the end of the introductory section in~\cite{kharlampovich-k-m:cayley}.  The theorem below answers a question raised implicitly at the end of~Section~8 in~\cite{kharlampovich-k-m:cayley}. 

\begin{theorem}\label{t:conjugacy}
There are Cayley automatic groups with undecidable Conjugacy Problem. 
\end{theorem}

Our last result concerns the Isomoprhism Problem. 

\begin{theorem}\label{t:ip}
The Isomorphism Problem is not decidable in the class of Cayley automatic groups. 
\end{theorem}

Our results follow from several results of Kharlampovich, Khoussainov, and Miasnikov~\cite{kharlampovich-k-m:cayley}, Bogopolski, Martino, and Ventura~\cite{bogopolski-m-v:cp}, and Levitt~\cite{levitt:ip}. 

Bogopolski, Martino, and Ventura proved that certain group extensions have decidable Conjugacy Problem (here and thereafter $F_n$ denotes the free group of rank $n$). 

\begin{theorem}[Corollary~7.6.~\cite{bogopolski-m-v:cp}]\label{t:bmv}
There exists a group of the form $\ZZ^d \rtimes_\tau F_n$ with undecidable Conjugacy Problem. 
\end{theorem}

The homomorphism $\tau: F_n \to \GL_d(\ZZ)$ constructed in the proof of Theorem~\ref{t:bmv} in~\cite{bogopolski-m-v:cp} is not injective. In fact, the image $\tau(F_n)$ is not finitely presented and the question of existence of a group of the form $\ZZ^d \rtimes_\tau F_n$ with undecidable Conjugacy Problem such that $\tau(F_n)$ is finitely presented was left open. A modification of the construction from~\cite{bogopolski-m-v:cp} that was provided in~\cite{sunic-v:unsolvable-cp} resolved this question.  

\begin{theorem}[\cite{sunic-v:unsolvable-cp}]\label{t:sv}
There exists a group of the form $\ZZ^d \rtimes_\tau F_n$ with undecidable Conjugacy Problem such that $\tau$ is injective.  
\end{theorem}

The primary goal of~\cite{sunic-v:unsolvable-cp} was to prove that the Conjugacy Problem is not decidable in the class of automaton groups (these are self-similar groups of rooted regular tree automorphisms generated by finite, invertible, synchronous transducers; see~\cite{grigorchuk-n-s:automata}). The class of automaton groups should not be confused with the class of automatic groups, as defined in~\cite{epstein-al:wp}, nor with its generalization, the Cayley automatic groups, as defined in~\cite{kharlampovich-k-m:cayley}. At present, the relation between the class of automaton groups and the class of Cayley automatic groups is not clear and only the latter is the subject of consideration in this work. 

The final ingredient in the proof of Theorem~\ref{t:conjugacy} is Theorem~\ref{t:semidirect}, which is not stated in~\cite{kharlampovich-k-m:cayley} in the form in which we quote it here, but it is a corollary of the other results presented there. Theorem~\ref{t:conjugacy} directly follows from Theorem~\ref{t:bmv} and Theorem~\ref{t:semidirect}. 

\begin{theorem}[\cite{kharlampovich-k-m:cayley}]\label{t:semidirect}
All groups of the form $\ZZ^d \rtimes F_n$ are Cayley automatic. 
\end{theorem}

As a direct corollary of Theorem~\ref{t:conjugacy} and the following result, we obtain Theorem~\ref{t:biautomatic}.  

\begin{theorem}[Theorem~8.5.~\cite{kharlampovich-k-m:cayley}]
Cayley biautomatic groups have decidable Conjugacy Problem. 
\end{theorem} 

As a direct corollary of Theorem~\ref{t:semidirect} and the following result of Levitt, we obtain Theorem~\ref{t:ip}. 

\begin{theorem}[\cite{levitt:ip}] 
The Isomorphism Problem is not decidable in the class of groups of the form $\ZZ^d \rtimes F_n$. 
\end{theorem}

It is important to observe that our examples of Cayley automatic groups that are not biautomatic and have undecidable Conjugacy Problem are not automatic in the standard sense. Indeed, one can prove the following. 

\begin{theorem}\label{t:subexponential}
If a group of the form $\ZZ^d \rtimes F_n$ has subexponential Dehn function, then it has decidable Conjugacy Problem. 
\end{theorem}

In the remaining sections we provide the necessary definitions and other details. 


\section{Cayley automatic and Cayley biautomatic groups}

Let $\SS$ be a finite alphabet. We will sometimes extend this alphabet by a special symbol $\diamond$ that is not in $\SS$, and we denote $\Sd = \SS \cup \{\diamond\}$. 

For an $n$-tuple of words $(w_1,\dots,w_n)$ over $\Sigma$ define the \emph{convolution} $\otimes(w_1,\dots,w_n)$ to be the word of length $\max\{|w_1|,\dots,|w_n|\}$ over $\Sdn$ in which the $j$-th symbol is $(\sigma_1,\dots,\sigma_n)$, where 
\[
\sigma_i = 
  \begin{cases} 
    \text{the } j \text{-th symbol of } w_i, & \text{ if } j \leq |w_i| \\
    \diamond, &  \text{otherwise}
  \end{cases}.
\]
For instance, 
\[
 \otimes(aaa,babaa,\emptyset) = 
 \begin{pmatrix} a \\ b \\ \diamond \end{pmatrix} 
 \begin{pmatrix} a \\ a \\ \diamond \end{pmatrix}
 \begin{pmatrix} a \\ b \\ \diamond \end{pmatrix}
 \begin{pmatrix} \diamond \\ a \\ \diamond \end{pmatrix}
 \begin{pmatrix} \diamond \\ a \\ \diamond \end{pmatrix}, 
\]
where $\emptyset$ denotes the empty word and the symbols in $\Sdn$ are written, for convenience, as columns. 

Let $R$ be an $n$-ary relation on $\SS^*$. The \emph{convolution} $\otimes R$ of $R$ is the language over $\Sdn$ defined by 
\[
 \otimes R = \{ \ \otimes(w_1,\dots,w_n) \mid (w_1,\dots,w_n) \in R \ \}. 
\] 
A relation $R$ is \emph{regular} over $\SS$ if its convolution $\otimes R$ is a regular language over $\Sdn$, i.e., $\otimes R$ is recognizable by a finite automaton over the alphabet $\Sdn$ (let us note that, in this work, the automata always read words from left to right). 

Let $G$ be a finitely generated group with finite generating set $S$. The right Cayley graph of $G$ with respect to $S$ is the graph $\Gamma(G,S)$ with $G$ as the set of vertices and, for each $g$ in $G$ and $s$ in $S$, an edge from $g$ to $gs$. The Cayley graph can be interpreted as a system of $|S|$ binary relations $E_s$ on $G$, for $s$ in $S$, where 
\[
 E_s = \{ \ (g,gs) \mid g \in G \ \}. 
\] 

A map $\bar{~}: G \to \SS^*$ induces $|S|$ binary relations on $\SS^*$ given by 
\[
 \ov{E}_s = \{ \ (\ov{g},\ov{gs}) \mid g \in G \ \}. 
\]

\begin{definition}
A finitely generated group $G$ with finite generating set $S$ is \emph{Cayley automatic} if there exists a finite alphabet $\SS$ and an injective map $\bar{~}:G \to \SS^*$ such that 

$\ov{G}$ is regular (over $\SS$) and

$\ov{E}_s$ is regular (over $\SS$), for every $s$ in $S$. 

In such a case the tuple $(\ov{G},\ov{E}_{s_1},\dots,\ov{E}_{s_k})$ is called an \emph{automatic structure} of the Cayley graph $\Gamma(G,S)$ or Cayley automatic structure of $G$ (with respect to $S=\{s_1,\dots,s_k\}$).  
\end{definition}

In addition to the right Cayley graph one may consider the left Cayley graph $\Gamma^\ell(G,S)$ as well. The vertex set is $G$ and, for each $g$ in $G$ and $s$ in $S$, an edge from $g$ to $sg$. The left Cayley graph can be interpreted as a system of $|S|$ binary relations $E^\ell_s$ on $G$, for $s$ in $S$, where 
\[
 E^\ell_s = \{ \ (g,sg) \mid g \in G \ \}. 
\] 

\begin{definition}
A finitely generated group $G$ with finite generating set $S$ is \emph{Cayley biautomatic} if there exists a finite alphabet $\SS$ and an injective map $\bar{~}:G \to \SS^*$ such that 

$\ov{G}$ is regular (over $\SS$), 

$\ov{E}_s$ is regular (over $\SS$), for every $s$ in $S$, and

$\ov{E}^\ell_s$ is regular (over $\SS$), for every $s$ in $S$. 

In such a case the tuple $(\ov{G},\ov{E}_{s_1},\dots,\ov{E}_{s_k},\ov{E}^\ell_{s_1},\dots,\ov{E}^\ell_{s_k})$ is called a \emph{biautomatic structure} of the pair of Cayley graphs $\Gamma(G,S)$ and $\Gamma^\ell(G,S)$ or Cayley biautomatic structure of $G$ (with respect to $S=\{s_1,\dots,s_k\}$).  
\end{definition}

It is important to observe that being Cayley automatic is property of the group and does not depend on the chosen finite generating set $S$, i.e., $G$ is Cayley automatic with respect to a finite generating $S$ if and only if it is Cayley automatic with respect to any of its other finite generating sets (Theorem~6.9.~\cite{kharlampovich-k-m:cayley}).  

All (bi)automatic groups, as defined in~\cite{epstein-al:wp} are Cayley (bi)automatic (Proposition~7.3. and Proposition~8.4.~\cite{kharlampovich-k-m:cayley}). The class of Cayley automatic groups is much wider than the class of automatic groups. For instance, it includes the Heisenberg group $H=\langle a,b \mid [a,[a,b]]=[b,[a,b]]=1\rangle$ and many other nilpotent groups that are not automatic (Example~6.6~\cite{kharlampovich-k-m:cayley}). Nevertheless, the class of Cayley (bi)automatic groups retains many algorithmic properties of (bi)automatic groups. For instance, every Cayley automatic group has Word Problem decidable in quadratic time and every Cayley biautomatic group has decidable Conjugacy Problem (Theorem~8.2. and Theorem~8.5.~\cite{kharlampovich-k-m:cayley}). 

The class of Cayley automatic groups has good closure properties. The following is, in particular, relevant for our purposes. 

\begin{theorem}[Theorem~10.3~\cite{kharlampovich-k-m:cayley}]\label{t:semidirect-general}
Let $A$ and $B$ be graph automatic groups with finite generating sets $X$ and $Y$, respectively. Let $\tau : B \to \Aut(A)$ be a homomorphism such that the automorphism $\tau(y)$ is automatic for every $y$ in $Y$. Then
the semidirect product $G = A \rtimes_\tau B$ is Cayley automatic.
\end{theorem}

By definition, an automorphism $\alpha$ of the Cayley automatic group $A$, which is automatic over $\SS$, is \emph{automatic} if the graph relation $\ov{\alpha}= \{(\ov{a},\ov{a^\alpha}) \mid a \in A\}$ induced by $\alpha$ is a regular relation over $\SS$. 

The semidirect product $A \rtimes_\tau B$ is the set of all pairs $(b,a)$, for $b \in B$, $a \in A$, with product defined by $(b_1,a_1)(b_2,a_2)=(b_1b_2,a_1^{\tau(b_2)}a_2)$. 

It is known that the free abelian group $A=\ZZ^d$ and the free group $B=F_n$ of finite ranks are automatic, and hence they are Cayley automatic. The argument in the proof of Proposition~10.5 in~\cite{kharlampovich-k-m:cayley}, showing that every automorphism of $\ZZ^2$ is automatic, can be extended to show that every automorphism of $\ZZ^d$ is automatic. In other words, multiplication of $d$-tuples of integers by any fixed $d \times d$ matrix in $\GL_d(\ZZ)$ is automatic. These observations, together with Theorem~\ref{t:semidirect-general} immediately imply Theorem~\ref{t:semidirect}.  


\section{(Free-abelian)-by-free groups with undecidable Conjugacy Problem}

For the duration of this section, let $A=\ZZ^d$ and $B=F_n$ (this agreement is not crucial for all statements, but this is the setting we are aiming for and there is no need to go into more general considerations). 

Let $C$ be a subgroup of $\Aut(A) = \GL_d(\ZZ)$. We say that $C$ has undecidable Orbit Problem if there is no algorithm that decides, on input consisting of arbitrary pair of vectors $u$ and $v$ in $A$, if there exists a matrix $c$ in $C$ such that $u^c=v$ (we use the right action of matrices on vectors; this is just the multiplication of vectors by the matrix $c$ on the right). Let $\tau: B \to \GL_d(\ZZ)$ be a homomorphism such that $\tau(B)=C$. If $C$ has undecidable Orbit Problem then $G = A \rtimes_\tau B$ has undecidable Conjugacy Problem. Indeed, as observed in~\cite{bogopolski-m-v:cp}, two vectors $u$ and $v$ in $A$ are conjugate in $G$ if and only if they are in the same orbit under the action of $C=\tau(B)$, and since the latter problem is undecidable, so is the Conjugacy Problem in $G$. 

A good way to construct orbit undecidable subgroups of $\GL_d(\ZZ)$ is provided in~\cite{bogopolski-m-v:cp} (Section~7; in particular, Proposition~7.5. and Corollary~7.6., the latter of which is listed in our introduction as Theorem~\ref{t:bmv}). Let $d \geq 4$ and let $H$ be a finitely presented group with undecidable Word Problem. Use the Mikhailova construction to obtain the corresponding finitely generated subgroup $H'$ of $F_2 \times F_2$ with undecidable Membership Problem and then consider $F_2 \times F_2$ as a subgroup of $\GL_d(\ZZ)$ through a specific embedding ($F_2 \times F_2$ embeds in $\GL_d(\ZZ)$, for $d \geq 4$) that turns the undecidability of the Word Problem in $H$ into undecidability of the Orbit Problem for $H' = C \leq \GL_d(\ZZ)$ (a specific embedding of $F_2 \times F_2$ with this property is spelled out precisely in~\cite{bogopolski-m-v:cp}). The group $G=A \rtimes_\tau B$, where $\tau : B \to \GL_d(\ZZ)$ is any homomorphism with $\tau(B)=C$, has undecidable Conjugacy Problem.   

The group $C$, as defined above, is finitely generated and not finitely presented. Thus, $\tau$ is not injective for any group of the form $G = A \rtimes_\tau B$ with $C=\tau(B)$ as in the above construction. 

The following modification, introduced in~\cite{sunic-v:unsolvable-cp}, provides groups of the form $G = A \rtimes_\tau B$ with undecidable Conjugacy Problem and injective $\tau$. Let $C=\langle g_1,\dots,g_n \rangle$ be an orbit undecidable subgroup of $\GL_d(\ZZ)$, let $B=F(f_1,\dots,f_n)$ be free of rank $n$, and let $C'=\langle g_1',\dots,g_n' \rangle$ be any free subgroup of rank $n$ of $\GL_2(\ZZ)$ (such subgroups do exist for any rank). Define $\tau: B \to \GL_{d+2}(\ZZ)$ by 
\[
 \tau(f_i) = \begin{bmatrix} g_i & 0_{d \times 2} \\ 0_{2 \times d} & g_i' \end{bmatrix}, 
\]  
for $i=1,\dots,n$, where $0_{d \times 2}$ and $0_{2 \times d}$ are the zero matrices of appropriate sizes. In other words, the action of $\tau(f_i)$ on the first $d$ coordinates of $\ZZ^{d+2}$ is the same as the action of the matrix $g_i$, and on the last two coordinates as the action of the matrix $g_i'$. It is clear that $\tau$ is injective (since it is injective ``on the last two coordinates''). Moreover, the undecidability of the Orbit Problem for $C$ in $\GL_d(\ZZ)$ induces the undecidabiluty of the Orbit Problem for the free subgroup $C'=\tau(B)$ in $\GL_{d+2}(\ZZ)$ (see Proposition~1 in~\cite{sunic-v:unsolvable-cp}, which is listed as Theorem~\ref{t:sv} in our introduction). 

At the end, we show that our examples of Cayley automatic groups that are not Cayley biautomatic and have undecidable Conjugacy Problem are not automatic under the standard definition. In fact, Theorem~\ref{t:subexponential} claims that they cannot even have subexponential Dehn functions. 

\begin{proof}[Proof of Theorem~\ref{t:subexponential}]
Let $G=\ZZ^d \rtimes_\tau F_n$ be a group with subexponential Dehn function. Bridson showed that the Dehn function of $G$ can be either polynomial or exponential and the former is possible only when $F_n$ has a subgroup $F'$ of finite index such that $\tau(F')$ is unipotent~\cite{bridson:z-by-f}. This implies that $\tau(F_n)$ is virtually solvable. Since virtually solvable subgroups of $\GL_d(\ZZ)$ have decidable Orbit Problem, it follows that $G$ has decidable Conjugacy Problem (see Proposition~6.9 and Corollary 6.10 in~\cite{bogopolski-m-v:cp}). 
\end{proof}

\newcommand{\etalchar}[1]{$^{#1}$}
\def\cprime{$'$}


\end{document}